\documentclass[runningheads]{llncs}
\usepackage[T1]{fontenc}
\usepackage{graphicx,verbatim}

\usepackage{booktabs}
\usepackage{array}
\usepackage[hidelinks,breaklinks=true]{hyperref}
\begin{document}

\title{3D Computational Modeling of the Vascular Heat-Sink Effect in Microwave Ablation of Tumors Close to the Aorta}
\titlerunning{Vascular Heat-Sink Effect in Microwave Ablation}

\begin{comment}  %% Removed for anonymized MICCAI submission
\author{First Author\inst{1}\orcidID{0000-1111-2222-3333} \and
Second Author\inst{2,3}\orcidID{1111-2222-3333-4444} \and
Third Author\inst{3}\orcidID{2222--3333-4444-5555}}
%
\authorrunning{F. Author et al.}
% First names are abbreviated in the running head.
% If there are more than two authors, 'et al.' is used.
%
\institute{Princeton University, Princeton NJ 08544, USA \and
Springer Heidelberg, Tiergartenstr. 17, 69121 Heidelberg, Germany
\email{lncs@springer.com}\\
\url{http://www.springer.com/gp/computer-science/lncs} \and
ABC Institute, Rupert-Karls-University Heidelberg, Heidelberg, Germany\\
\email{\{abc,lncs\}@uni-heidelberg.de}}

\end{comment}

% \author{Anonymized Authors}  %% Added for anonymized MICCAI submission
% \authorrunning{Anonymized Author et al.}
% \institute{Anonymized Affiliations \\
%     \email{email@anonymized.com}}

\author{Youjin Kim\inst{1} \and
Kyungho Yoon\inst{2}\and Minwoo Shin\inst{1}\thanks{Corresponding author.}}
\authorrunning{Y. Kim et al.}
\institute{Department of Software, Yonsei University (Mirae Campus), Wonju, 26493, Republic of Korea\and School of Mathematics and Computing, Yonsei University, Seoul, 03722, Republic of Korea\\
\email{mshin@yonsei.ac.kr}}    
  
\maketitle              % typeset the header of the contribution
\begin{abstract}
Microwave ablation (MWA) is a percutaneous ablation that induces thermal necrosis by generating localized high-temperature regions within tumor tissue. However, when a tumor is located close to a large vessel such as the aorta, blood flow continuously removes heat from the ablation zone. This phenomenon causes the vascular heat-sink effect and thereby increases the risk of incomplete ablation. This study constructs an ANSYS-based numerical framework to investigate the heat-sink effect during MWA of tumors located close to the aorta. The computational model consists of the tumor, surrounding tissue, vessel wall, blood flow domain, and microwave antenna. Microwave-induced heating was represented by an equivalent heat source mapped from the High-Frequency Structure Simulator (HFSS)-derived volume loss density, rather than by direct electromagnetic-thermal coupling within the thermal-fluid solver. Two parametric studies were performed to quantify the heat-sink effect. First, the tumor–aorta distance was varied at \(D = 0~\mathrm{mm}\), \(4~\mathrm{mm}\), and \(8.5~\mathrm{mm}\) to assess how vessel proximity affects temperature elevation and ablation zone formation. 
Second, at each tumor-aorta distance, three blood-flow conditions were compared: no blood, stagnant blood ($V = 0$ m/s), and normal aortic blood flow ($V=0.65$ m/s). This comparison enabled the heat-sink effect of blood flow to be independently evaluated. The results highlight the clinical need for simulation-based treatment planning in patients with tumors located close to major vessels.

\keywords{Microwave ablation \and Heat-sink effect \and Aorta \and Blood flow \and Thermal-fluid simulation} % Authors must provide keywords and are not allowed to remove this Keyword section. 

\end{abstract}
\section{Introduction}

Microwave ablation (MWA) is a percutaneous thermal ablation technique that uses microwave energy to form a localized high-temperature region within tumor tissue, which causes irreversible damage and cell death through protein denaturation and irreversible cellular injury. During the procedure, a microwave antenna is inserted into the tumor, and electromagnetic energy is delivered to heat the surrounding tissue generally to a cytotoxic temperature of \(60^{\circ}\mathrm{C}\) or higher~\cite{cavagnaro2011antenna,wei2022spherical,carrafiello2008principles,lopresto2017properties}. MWA is a minimally invasive treatment with fewer complications and faster recovery than surgical resection, and it has been applied to various solid tumors such as hepatocellular carcinoma (HCC), lung cancer, and renal cell carcinoma (RCC)~\cite{carrafiello2008principles,deshazer2017experimental}.

For thermal ablation therapy to be effective, the entire tumor and a sufficient safety margin must reach the target ablation temperature. However, when a tumor is adjacent to a large vessel such as the aorta, continuous blood flow removes heat from the surrounding tissue and suppresses the temperature increase near the vessel wall~\cite{keangin2013heat,ringe2015heatsink,yu2008veinsize,siriwardana2017perfusion}. This vascular heat-sink effect can cause incomplete ablation of the tumor, which may be a major cause of local tumor recurrence after treatment~\cite{keangin2013heat,ringe2015heatsink,siriwardana2017perfusion}.

The thermal distribution formed during MWA is affected by various factors such as the tumor-to-vessel distance, tumor size and morphology, and blood flow rate. These variables vary substantially among patients, and it is difficult to sufficiently predict heat transfer behavior and cooling effects using real-time imaging alone. Therefore, simulation-based treatment planning using computational analysis tools is useful for predicting the thermal distribution, cooling near the vessel, and expected ablation range before the procedure~\cite{keangin2013heat,gas2012temperature,prakash2010hepatic,heshmat2024patient}. In particular, a coupled thermal-fluid model enables quantitative evaluation of whether the entire tumor reaches the target temperature, the degree of cooling in the vessel-adjacent region, and the final extent of ablation zone formation.

The objective of this study was to reproduce the heat-sink effect occurring during MWA of a tumor adjacent to the aorta using ANSYS-based coupled solid-fluid numerical simulations and to analyze changes in tumor temperature distribution and ablation region according to vessel proximity and blood flow conditions. For this purpose, the tumor-to-aorta distance was set to D = 0 mm, 4 mm, and 8.5 mm, and no blood, stagnant blood ($V = 0$ m/s), and normal aortic blood flow ($V = 0.65$ m/s) conditions were applied. Then, the contour figures and ablation zone volumes of each case were compared to identify the risk of incomplete ablation in tumors adjacent to large vessels and to present the usefulness of simulation-based treatment planning. 

\begin{figure}[t]
    \centering

    \begin{minipage}[c]{0.431\linewidth}
        \centering
        \includegraphics[
            width=\linewidth,
            angle=90,
            origin=c
        ]{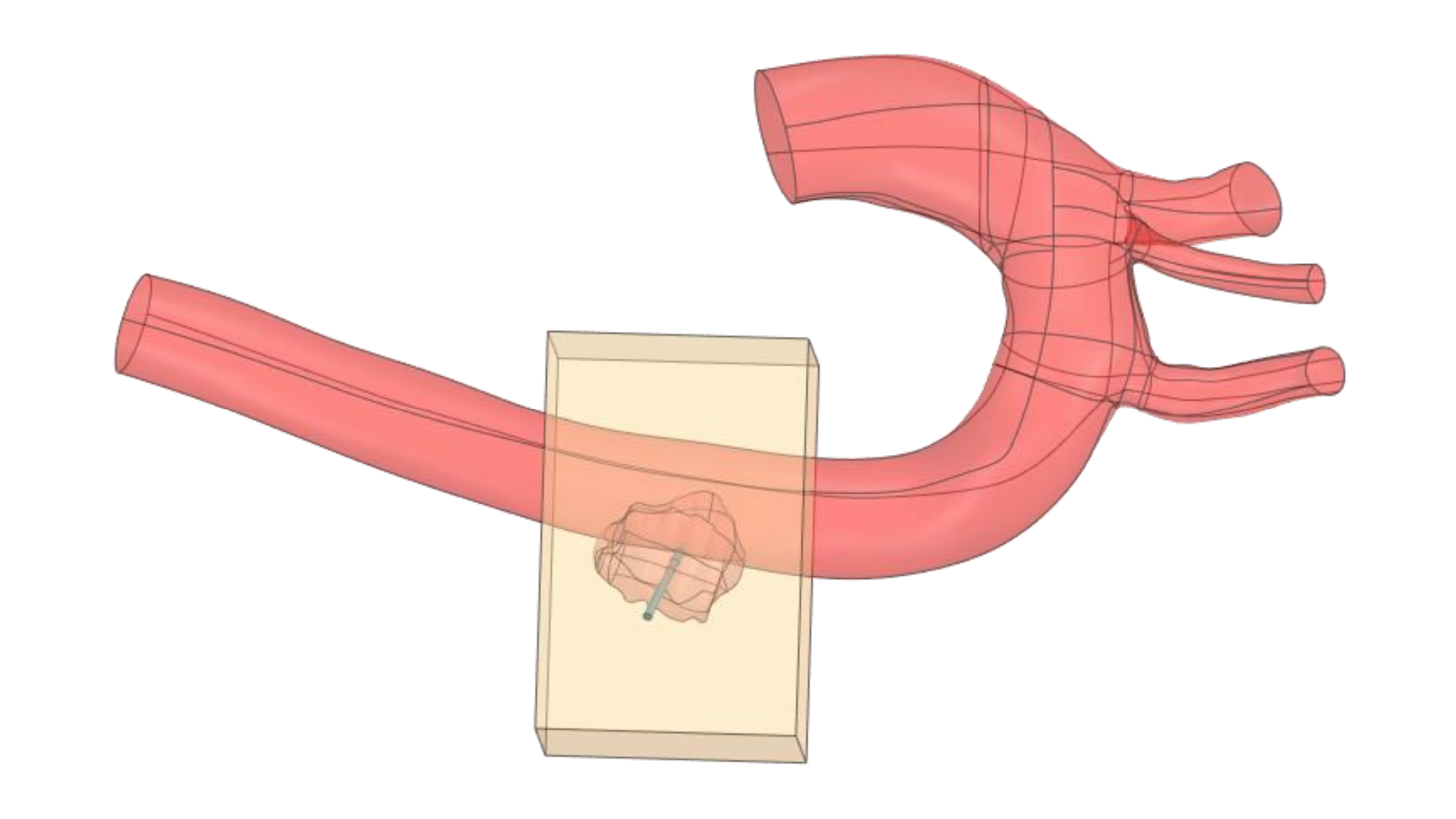}

        {\small (a) Geometry model}
        \label{fig:geometry_model_original}
    \end{minipage}
    \hfill
    \begin{minipage}[c]{0.37\linewidth}
        \centering
        \includegraphics[
            width=\linewidth,
            origin=c
        ]{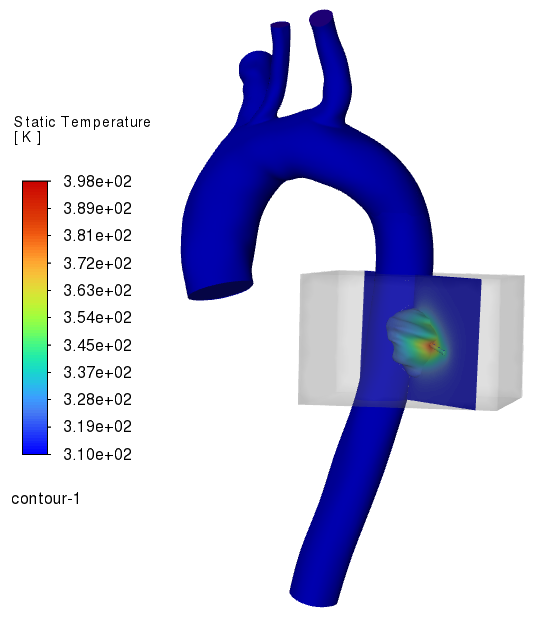}

        {\small (b) Temperature distribution}
        \label{fig:geometry_model_calculated}
    \end{minipage}

    \caption{Geometry models used for the MWA simulation.
    (a) Overall geometry model including the tumor, surrounding tissue, aorta, and microwave antenna.
    (b) Calculated geometry model showing the temperature distribution after \(300~\mathrm{s}\) of \(20~\mathrm{W}\) microwave heating under the contact condition \((D = 0~\mathrm{mm})\) with normal aortic blood flow \((V = 0.65~\mathrm{m/s})\).}
    \label{fig:geometry_models}
\end{figure}

\section{Method}

\begin{figure}[t]
    \centering
    \includegraphics[width=0.5\linewidth]{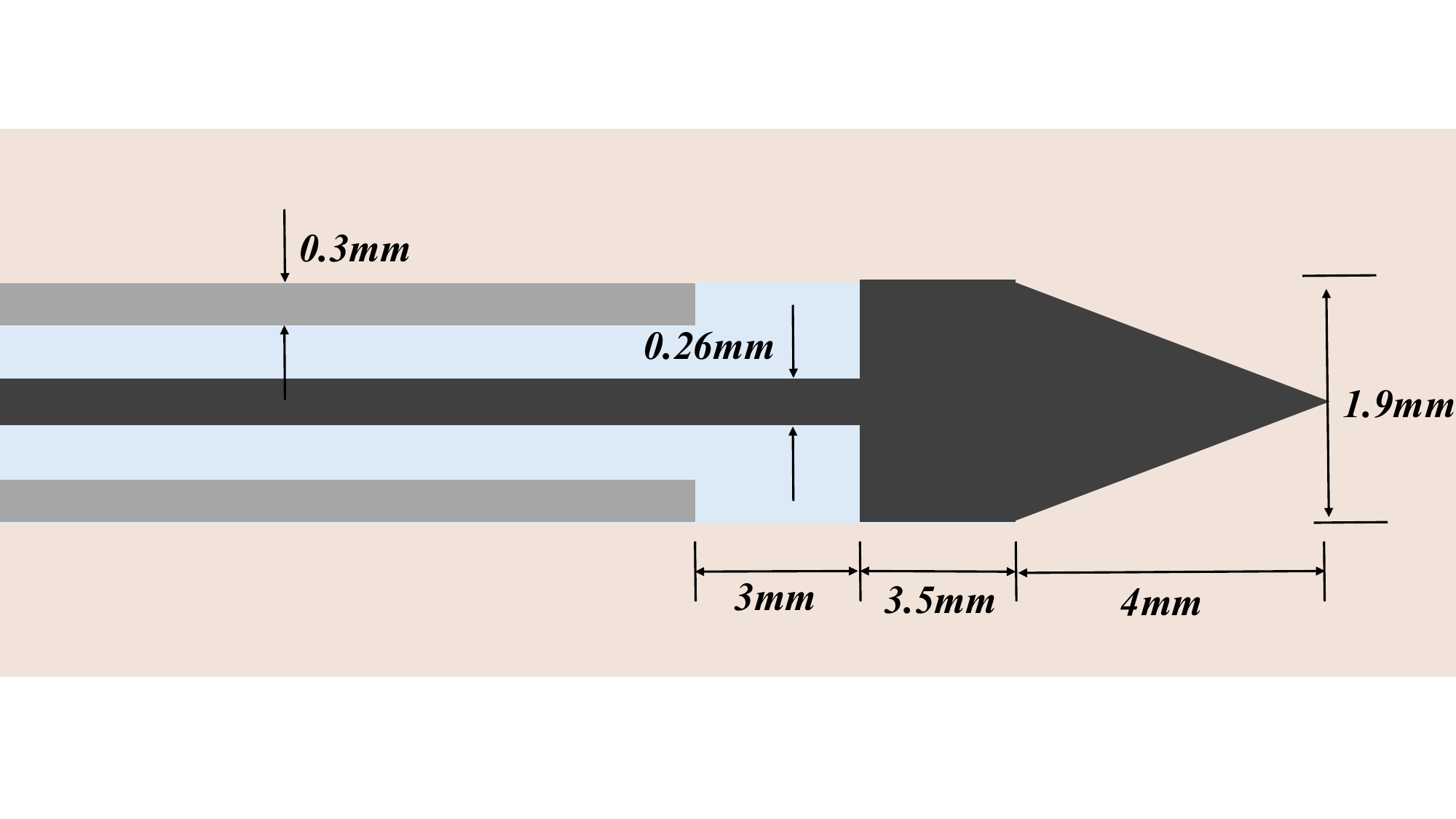}
    \caption{Configuration of the microwave antenna with an asymmetric dipole structure.}
    \label{fig:antenna_geometry}
\end{figure}

\subsection{Geometry Construction}

The computational model in this study consisted of a tumor, surrounding tissue, aorta, blood fluid domain, and microwave antenna. The overall geometry model and the calculated geometry model are shown in Fig.~\ref{fig:geometry_models}. The tumor geometry was obtained from the 3D-IRCADb-01 liver segmentation dataset provided by IRCAD France\footnote{\url{https://www.ircad.fr/research-and-development/data-sets/liver-segmentation-3d-ircadb-01/}}. A scale factor of 0.035 was applied to reproduce the tumor size characteristics reported in clinical MWA literature on patients with retroperitoneal tumors. The reported tumor size had a mean diameter of 2.6 cm and a range of \(0.5\)--\(4.5~\mathrm{cm}\). The scaled tumor model was placed within the surrounding tissue adjacent to the aorta to represent the tumor tissue treated by MWA.

The aorta geometry was obtained as an STL-based vascular model from the Vascular Model Repository\footnote{\url{https://www.vascularmodel.com/index.html}}. The original model was constructed from data from a 26-year-old healthy male. In this study, the downloaded vascular geometry was used as the final computational domain after a scale factor of 0.472 was applied. This scaling was determined to correspond to the normal thoracic aortic diameter range reported by CMR reference values, including an ascending aorta diameter of \(27 \pm 4~\mathrm{mm}\) and a distal descending aorta diameter of \(18 \pm 3~\mathrm{mm}\) in adult males.

The microwave antenna was modeled as a slot-type coaxial MWA applicator based on the KY-2450B configuration and a tapered-tip antenna geometry reported in previous MWA antenna studies~\cite{wei2022spherical,gas2012temperature}. The antenna configuration is shown in Fig.~\ref{fig:antenna_geometry}. The antenna consisted of an inner conductor, Teflon dielectric layer, outer conductor, exposed Teflon slot, distal metallic cylinder, and tapered metallic tip. The inner conductor was aligned with the antenna axis and acted as the microwave feed conductor. The surrounding Teflon layer electrically insulated the inner conductor from the outer conductor. The outer conductor enclosed the coaxial shaft and ended before the distal end. This structure exposed a 3-mm Teflon slot, which served as the main microwave emission region.

Distal to the exposed Teflon slot, a 3.5-mm metallic cylinder and a 4-mm tapered conical tip were constructed to reproduce the distal applicator geometry. The maximum tip diameter was 1.9 mm. In the antenna cross-section, the outer conductor thickness was 0.3 mm, and the dielectric spacing around the inner conductor was 0.26 mm. Copper was assigned to the inner conductor and distal metallic tip, stainless steel to the outer conductor, and Teflon to the dielectric regions.

When blood flow conditions were applied, the internal region of the aorta was defined as the blood fluid domain. Under the no blood condition, the blood fluid domain was excluded to simulate a state in which cooling by blood flow did not occur. In addition, tumor-to-vessel distances of \(D = 0~\mathrm{mm}\), \(4~\mathrm{mm}\), and \(8.5~\mathrm{mm}\) were defined as the main spatial variables. This model design allowed quantitative comparison of the effects of tumor--vessel proximity on the temperature field and ablation zone formation during MWA.

\begin{table}[t]\footnotesize
\centering
\caption{Material properties used in the computational model. Here, $\rho$ is the density, $C_p$ is the specific heat, $k$ is the thermal conductivity, and $\mu$ is the viscosity.}
\label{tab:material_properties}
\begin{tabular}{lcccc}
\toprule
\textbf{Material} 
& \textbf{$\rho$} 
& \textbf{$C_p$} 
& \textbf{$k$} 
& \textbf{$\mu$} \\
& \textbf{[kg/m$^3$]} 
& \textbf{[J/(kg$\cdot$K)]} 
& \textbf{[W/(m$\cdot$K)]} 
& \textbf{[Pa$\cdot$s]} \\
\midrule
Blood           & 1058 & 3960 & 0.45  & 0.0035 \\
Tissue          & 1030 & 3600 & 0.497 & - \\
Tumor           & 1040 & 3960 & 0.57  & - \\
Teflon          & 2200 & 1300 & 0.25  & - \\
Copper          & 8930 & 385  & 401   & - \\
Stainless steel & 7920 & 502  & 15    & - \\
\bottomrule
\end{tabular}
\end{table}

\subsection{Material Properties} 
The thermal and fluid properties assigned to these materials are listed in
Table~\ref{tab:material_properties}. The biological material properties of
tissue, tumor, and blood were selected based on previous MWA simulations
and reported thermal-property ranges for biological tissues~\cite{radjenovic2021efficacy}. Copper, stainless steel,
and Teflon were assigned to the antenna components following material
selections used in previous MWA probe simulations~\cite{gorman2022numerical},
while their numerical thermal properties were assigned from the ANSYS
material database.

\subsection{HFSS-to-Fluent Heat Source Mapping}

Microwave heating was applied by transferring the volume loss density from HFSS to the energy source term in Fluent. The volume loss density at $2.45~\mathrm{GHz}$ and $20~\mathrm{W}$ was exported from HFSS. Python was used to convert the data into a structured Cartesian grid. A UDF in Fluent read the grid data and found the nearest volume loss density value for each cell centroid. The source term was set to zero for cells outside the grid range. A temperature-dependent attenuation factor was also applied. The heat source remained unchanged below $333.15~\mathrm{K}$. It decreased linearly from $1.0$ to $0.50$ between $333.15~\mathrm{K}$ and $373.15~\mathrm{K}$. It remained at $0.50$ above $373.15~\mathrm{K}$. The effective heat source was calculated by multiplying the raw volume loss density by the attenuation factor. This one-way HFSS-to-Fluent mapping applied a non-uniform heat source near the antenna slot. The electromagnetic field was not recalculated during the transient thermal-fluid analysis.

\subsection{Physical Model} 
A coupled thermal-fluid simulation was performed to reproduce the vascular heat-sink effect in tumors adjacent to the aorta. In solid regions such as the tissue, tumor, and antenna, temperature changes due to heat conduction were calculated. Although previous MWA studies have reported that perfusion and vessel-adjacent blood flow can reduce lesion size or deform thermal lesions near vessels~\cite{ringe2015heatsink,yu2008veinsize,siriwardana2017perfusion}, they did not explicitly isolate the effects of tumor--aorta distance and blood flow condition in an aorta-adjacent tumor model. In this study, convective heat removal was therefore modeled in the blood fluid domain to quantitatively evaluate the influence of large-vessel cooling on tumor temperature distribution and ablation zone formation. For the flow analysis of blood, the SST \(k\)-\(\omega\) model was applied to account for turbulence effects under aortic blood flow conditions.

\section{Experiments}
\subsection{Tumor-to-Aorta Distance Study}

To analyze the effect of tumor--vessel proximity on the vascular heat-sink effect, tumor-to-aorta distances of \(D = 0~\mathrm{mm}\), \(4~\mathrm{mm}\), and \(8.5~\mathrm{mm}\) were defined with reference to previous MWA heat-transfer and heat-sink studies that evaluated temperature variation and ablation-zone deformation near vessels at different spatial positions~\cite{keangin2013heat,ringe2015heatsink,yu2008veinsize}. The \(D = 0~\mathrm{mm}\) condition was used as the direct-contact case between the tumor and the aorta. The \(D = 4~\mathrm{mm}\) and \(D = 8.5~\mathrm{mm}\) conditions were defined as cases with increased tissue gaps between the tumor and the aorta. Under each distance condition, tumor temperature distribution, vessel-side temperature reduction, and ablation zone volume were compared to evaluate whether the vascular heat-sink effect became stronger as the tumor was located closer to the aorta. In particular, the formation of an asymmetric ablation zone and the reduction of the \(60^{\circ}\mathrm{C}\) isothermal region near the vessel side were analyzed to assess the effect of tumor-to-aorta distance on incomplete ablation risk.

\subsection{Blood Flow Condition Study} 
To compare the effect of blood-flow-induced cooling under different blood flow conditions, no blood, stagnant blood, and normal aortic blood flow conditions were defined. The no blood condition was used as the reference case without a blood fluid domain, in which convective heat removal by blood flow did not occur. The stagnant blood condition was defined as the case in which the blood domain was present but the flow velocity was set to \(V = 0~\mathrm{m/s}\). The normal aortic blood flow condition was defined with \(V = 0.65~\mathrm{m/s}\), which approximates the reported aortic blood-flow velocity measured by 4D-flow CMR~\cite{kawelboehm2020cmr,garcia2018aorta}. Under each blood flow condition, tumor temperature distribution, vessel-side temperature reduction, and ablation zone volume were compared to evaluate whether heat removal increased and the tumor temperature near the vessel decreased as blood flow velocity increased.

\begin{figure}[t]
    \centering
    \includegraphics[
        width=0.95\linewidth
    ]{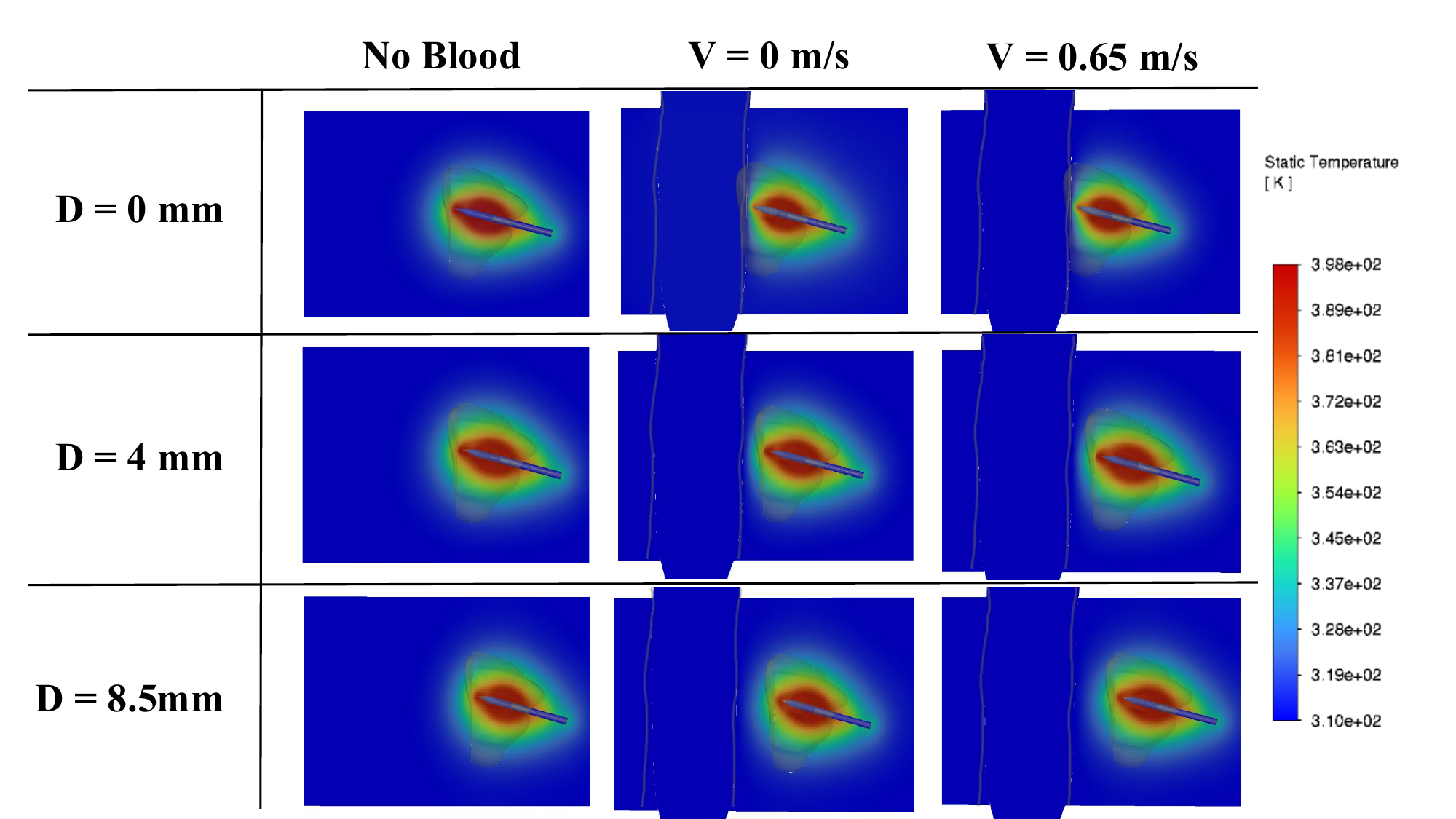}
    \caption{Temperature contours under different tumor--aorta distances and blood flow conditions.
    Each row represents a different tumor--aorta distance: \(D = 0~\mathrm{mm}\), \(4~\mathrm{mm}\), and \(8.5~\mathrm{mm}\).
    Each column represents a different blood flow condition: no blood, stagnant blood \((V = 0~\mathrm{m/s})\), and normal aortic flow \((V = 0.65~\mathrm{m/s})\).
    For the blood-present cases, the gray aortic outline is superimposed on the temperature contour to indicate the location of the aorta and to visualize the tumor--vessel separation distance.}
    \label{fig:contour_comparison}
\end{figure}

\begin{table}[t]\footnotesize
\centering
\caption{Quantitative comparison of ablation volumes under different tumor--aorta distances and blood flow conditions.}
\label{tab:ablation_volume_comparison}
\begin{tabular}{@{}c l c c c c@{}}
\hline
\(D\) (mm)
& Blood flow
& \(V_{60}^{\mathrm{total}}\) \((\mathrm{m^3})\)
& \(V_{100}^{\mathrm{total}}\) \((\mathrm{m^3})\)
& \(V_{60}^{\mathrm{tumor}}\) \((\mathrm{m^3})\)
& \(C_{60}^{\mathrm{tumor}}\) (\%) \\
\hline

\(0.0\) & No blood
& \(8.3149{\times}10^{-6}\)
& \(9.7897{\times}10^{-7}\)
& \(5.9341{\times}10^{-6}\)
& 66.44 \\

\(0.0\) & \(V = 0.00~\mathrm{m/s}\)
& \(6.5699{\times}10^{-6}\)
& \(7.8126{\times}10^{-7}\)
& \(4.6968{\times}10^{-6}\)
& 52.59 \\

\(0.0\) & \(V = 0.65~\mathrm{m/s}\)
& \(6.5748{\times}10^{-6}\)
& \(7.8166{\times}10^{-7}\)
& \(4.4392{\times}10^{-6}\)
& 49.71 \\

\hline

\(4.0\) & No blood
& \(8.3231{\times}10^{-6}\)
& \(9.8751{\times}10^{-7}\)
& \(5.9372{\times}10^{-6}\)
& 66.48 \\

\(4.0\) & \(V = 0.00~\mathrm{m/s}\)
& \(8.0132{\times}10^{-6}\)
& \(9.7711{\times}10^{-7}\)
& \(5.8269{\times}10^{-6}\)
& 65.25 \\

\(4.0\) & \(V = 0.65~\mathrm{m/s}\)
& \(8.0129{\times}10^{-6}\)
& \(9.7697{\times}10^{-7}\)
& \(5.8287{\times}10^{-6}\)
& 65.27 \\

\hline

\(8.5\) & No blood
& \(8.3839{\times}10^{-6}\)
& \(9.9546{\times}10^{-7}\)
& \(5.9740{\times}10^{-6}\)
& 66.89 \\

\(8.5\) & \(V = 0.00~\mathrm{m/s}\)
& \(8.3678{\times}10^{-6}\)
& \(1.0031{\times}10^{-6}\)
& \(5.9874{\times}10^{-6}\)
& 67.05 \\

\(8.5\) & \(V = 0.65~\mathrm{m/s}\)
& \(8.3696{\times}10^{-6}\)
& \(1.0030{\times}10^{-6}\)
& \(5.9766{\times}10^{-6}\)
& 66.93 \\

\hline
\end{tabular}
\end{table}

\section{Results}

\subsection{Distance-Dependent Comparison}
The temperature distributions under different tumor--aorta distances and blood flow conditions are shown in Fig.~\ref{fig:contour_comparison}. The quantitative ablation volumes are listed in Table~\ref{tab:ablation_volume_comparison}. In all cases, the highest temperature region appeared near the microwave antenna slot. Heat spread along the antenna axis into the tumor and surrounding tissue. However, the vessel-side temperature field changed with the tumor--aorta distance.

In the \(D = 0~\mathrm{mm}\) condition, the tumor was in direct contact
with the aorta in the blood-present models. In the no blood condition,
the high-temperature region expanded broadly around the antenna, and
\(C_{60}^{\mathrm{tumor}}\) reached 66.44\%. When the blood domain was
included, the high-temperature contour was compressed near the vessel
side. The \(C_{60}^{\mathrm{tumor}}\) value decreased to 52.59\% under
stagnant blood and 49.71\% under normal flow. These values correspond
to reductions of 13.85 and 16.73 percentage points relative to the no
blood condition. The residual tumor volume below
\(60^{\circ}\mathrm{C}\) increased from \(2997.02~\mathrm{mm^3}\) to
\(4233.57~\mathrm{mm^3}\) and \(4491.17~\mathrm{mm^3}\), respectively.
These results show that direct contact with the blood domain strongly
reduced vessel-side thermal expansion.

In the \(D = 4~\mathrm{mm}\) condition, a small tissue gap was present
between the tumor and the aorta. The no blood condition showed a broad
ablation region, with a \(C_{60}^{\mathrm{tumor}}\) value of 66.48\%.
The value decreased to 65.25\% under stagnant blood and remained at
65.27\% under normal aortic blood flow. The two blood-present values
differed by only 0.02 percentage points. The residual tumor volume below
\(60^{\circ}\mathrm{C}\) increased from \(2993.93~\mathrm{mm^3}\) in
the no blood condition to \(3103.40~\mathrm{mm^3}\) under stagnant blood
and \(3101.67~\mathrm{mm^3}\) under normal flow. Thus, the blood-present
models produced a small reduction in tumor coverage at
\(D = 4~\mathrm{mm}\), while the increase in blood velocity did not
produce a further reduction.

In the \(D = 8.5~\mathrm{mm}\) condition, the effect of blood flow on
the main ablation zone was limited. The central high-temperature region
was maintained near the antenna tip under all blood flow conditions.
As shown in Table~\ref{tab:ablation_volume_comparison},
\(C_{60}^{\mathrm{tumor}}\) remained between 66.89\% and 67.05\%.
The maximum difference was 0.16 percentage points. This result shows
that the vascular heat-sink effect on tumor coverage became weak when
the tumor--aorta distance increased to \(8.5~\mathrm{mm}\).

\subsection{Blood Flow Condition Comparison}

The same results were also compared by blood flow condition, as summarized in Table~\ref{tab:ablation_volume_comparison}. The no blood condition was used as the reference case. Since the blood fluid domain was excluded, convective heat removal by blood flow did not occur. Therefore, the ablation zone remained broad and relatively symmetric across all distance conditions.

The stagnant blood condition separated the effect of the blood domain from the effect of blood velocity. In this condition, the blood domain was present but the flow velocity was zero. The largest change appeared in the direct-contact case. This means that the presence of blood near the tumor was sufficient to absorb heat and reduce thermal expansion at the vessel side.

The normal aortic blood flow condition added the effect of blood
velocity. The additional reduction in tumor coverage was clearest in
the direct-contact condition. At \(D = 0~\mathrm{mm}\),
\(C_{60}^{\mathrm{tumor}}\) decreased from 52.59\% under stagnant blood
to 49.71\% under normal flow. At \(D = 4~\mathrm{mm}\) and
\(D = 8.5~\mathrm{mm}\), the differences between stagnant and normal
flow were only 0.02 and 0.12 percentage points, respectively. Thus, the
added effect of blood velocity on tumor coverage was small when a tissue
gap was present.

These results show that tumor--aorta distance and blood flow condition had different roles. Tumor--aorta distance controlled whether vessel cooling reached the tumor. Blood flow condition controlled how strongly heat was removed once the vessel-side thermal field interacted with the blood domain. Therefore, MWA treatment planning for tumors near large vessels should consider tumor coverage, the vessel-side \(60^{\circ}\mathrm{C}\) margin, and patient-specific vessel position in addition to input power and central temperature.

\section{Conclusion}
In this study, an ANSYS-based coupled thermal-fluid simulation was used to analyze the vascular heat-sink effect during MWA of tumors close to the aorta. Microwave heating was modeled as an equivalent localized heat source. The results showed that a shorter tumor-to-aorta distance and the presence of the blood domain reduced the temperature near the tumor-vessel interface and decreased tumor coverage. This effect was most clearly observed in the \(D = 0~\mathrm{mm}\) condition. A smaller
reduction was observed at \(D = 4~\mathrm{mm}\), while tumor coverage remained nearly unchanged at \(D = 8.5~\mathrm{mm}\). These findings show that tumor-to-aorta distance and blood flow should be considered in MWA treatment planning. They also suggest that simulation can help predict vessel-side heat loss before the procedure. Future work could incorporate varied or randomly generated tumor geometries and machine-learning techniques such as surrogate modeling to support patient-specific ablation planning.

\begin{credits}
\subsubsection{\ackname} 
This work was supported by the Korea Medical Device Development Foundation grant funded by the Korean government (the Ministry of Science and ICT, the Ministry of Trade, Industry and Energy, the Ministry of Health and Welfare, and the Ministry of Food and Drug Safety) (Grant No. RS-2026-25543484).

This research was also supported by the ANCHOR Program through the Gangwon ANCHOR Center, funded by the Ministry of Education (MOE) and Gangwon State (G.S.), Republic of Korea (Grant No. 2026-ANCHOR-10-006).

This research was further supported by the Ministry of Science and ICT (MSIT), Korea, under the National Program in Medical AI Semiconductor (Grant No. 2024-0-00096), supervised by the Institute of Information \& Communications Technology Planning \& Evaluation (IITP) in 2026.

\subsubsection{\discintname}
The authors have no competing interests to declare that are relevant to the content of this article.
\end{credits}
%
% ---- Bibliography ----
%
% BibTeX users should specify bibliography style 'splncs04'.
% References will then be sorted and formatted in the correct style.
%
\bibliographystyle{splncs04}
\bibliography{references}

\end{document}